\documentclass[amssymb,11pt,reqno]{amsart}

\newcommand{\R}{\mathbb{R}}
\newcommand{\C}{\mathbb{C}}

\newcommand{\N}{\mathbb{N}}
\newcommand{\G}{\Gamma}
\newcommand{\rth}{R_{\theta}}
\newcommand{\test}{\mathcal S}

\newtheorem*{main}{Main Result}
\newtheorem{thm}{Theorem}

\newtheorem{lemma}{Lemma}
\newtheorem{prop}{Proposition}

\usepackage{graphicx}
\usepackage{latexsym}
\def\vol{\mbox{\rm Vol}}
\def\Re{\mbox{\rm Re}}

\begin{document}

\vspace{1.5 cm}

%\begin{center}{ Last modified Dec 2, 2004  } \end{center}

%%%%%%%%%%%%%%%%%%%%%%%%%%%%%%%%%%%%%%%%%%%%%%%%%%%%%%%%%%%%%%%%%%%%
\title[The  Modified complex Busemann-Petty problem on sections of convex bodies]
      {The  Modified complex Busemann-Petty problem on sections of convex bodies.}
      
\author{Marisa Zymonopoulou}
\address{
Marisa Zymonopoulou\\
Department of Mathematics\\
University of Missouri\\
Columbia, MO 65211, USA}

\email{marisa@@math.missouri.edu}
%%%%%%%%%%%%%%%%%%%%%%%%%%%%%%%%%%%%%%%%%%%%%%%%%%%%%%%%%%%%%%%%%%%%\begin{abstract}
\begin{abstract}
The complex Busemann-Petty problem asks whether origin symmetric convex bodies in $\C^n$ with smaller central hyperplane sections necessarily have smaller volume. The answer is affirmative if $n\leq 3$ and negative if $n \geq 4.$ Since the answer is negative in most dimensions, it is natural to ask what conditions on the $(n-1)$-dimensional volumes of the central sections of complex convex bodies with complex hyperplanes allow to compare the $n$-dimensional volumes. In this article we give necessary conditions on the section function in order to obtain an affirmative answer in all dimensions. The result is the complex analogue of [KYY].

\end{abstract}

\maketitle

\section{Introduction}

The Busemann-Petty problem was completely solved in the late 90's as a result of a series of papers of many mathematicians ([LR], [Ba], [Gi], [Bo], [Lu], [Pa], [Ga], [Zh1], [K1], [K2], [Zh2], [GKS]; see [K5, p.3] for the history of the solution).  The problem asks the following:

\noindent
Suppose $K$ and $L$ are two origin symmetric convex bodies in $\R^n$ such that for every $\xi \in S^{n-1},$

$$\vol_{n-1}\bigl(K\cap \xi^{\perp}\bigr) \leq \vol_{n-1}\bigl(L\cap \xi^{\perp}\bigr).$$
Does it follow that
$$\vol_{n}\bigl(K\bigr) \leq \vol_{n}\bigl(L \bigr) \ \ ?$$

The problem has an affirmative answer only if $n \leq 4.$ Since the answer is negative in most dimensions, it is natural to ask what conditions on the $(n-1)$-dimensional volumes of central sections do allow to compare the $n$-dimensional volumes.
Such conditions were found in [KYY]. The result is as
follows.

\noindent
For an origin symmetric convex body $K$ in $\R^n$ define the section function $$S_K(\xi)=\vol_{n-1}(K\cap \xi^{\perp}), \ \xi \in S^{n-1}.$$
Suppose $K$ and $L$ are origin symmetric convex smooth bodies  in $\R^n$ and $\alpha \in \R$ with $ \alpha \geq n-4.$ Then, the inequality
$$\bigl(-\Delta \bigr)^{\alpha/2}S_{K}(\xi)\leq \bigl(-\Delta \bigr)^{\alpha/2}S_{L}(\xi), \ \xi \in S^{n-1}$$
implies that $\vol_{n}(K)\leq \vol_{n}(L).$ If $\alpha <n-4$ this is not necessarily true. Here, $\Delta$ is the Laplace operator on $\R^n.$

In this article we study the complex version of this problem.
 For $\xi \in \C^n, \ |\xi|=1$ we denote by
$$H_{\xi}=\{z\in \C^n \ : \ (z,\xi)=\sum\limits_{k=1}^{n}z_k \overline{\xi}_k=0\}$$

\noindent
the complex hyperplane perpendicular to $\xi.$

Origin symmetric convex bodies in $\C^n$ are the unit balls of norms on $\C^n.$ We denote by $\|\cdot\|_K$ the norm corresponding to the body $K$
$$K=\{z \in \C^n \ : \|z\|_{K}\leq 1\}.$$

\noindent
We identify $\C^n$ with $\R^{2n}$ using the mapping
$$\xi=(\xi_1,\ldots ,\xi_n)=(\xi_{11}+i\xi_{12},\ldots ,\xi_{n1}+i\xi_{n2})\longmapsto (\xi_{11},\xi_{12},\ldots, \xi_{n1},\xi_{n2})$$
\noindent
and observe that under this mapping the complex hyperplane $H_{\xi}$ turns into a $(2n-2)$-dimensional subspace of $\R^{2n}$ orthogonal to the vectors

   \begin{center}
   $\xi=(\xi_{11},\xi_{12},\ldots, \xi_{n1},\xi_{n2})$   and   $\xi^{\perp}=(-\xi_{12},\xi_{11},\ldots, -\xi_{n2},\xi_{n1}).$
   \end{center}

\noindent
Since norms on $\C^n$ satisfy the equality
$$\|\lambda z\|=|\lambda|\|z\|, \ \ \forall z \in \C^n, \ \forall \lambda \in \C^n,$$
origin symmetric complex convex bodies correspond to those origin symmetric convex bodies $K$ in $\R^{2n}$ that are invariant with respect to any coordinate-wise two-dimensional rotation, namely for each $\theta \in [0,2\pi]$ and each $x=(x_{11},x_{12},\ldots, x_{n1},x_{n2}) \in \R^{2n}$
\begin{equation}
\|x\|_{K}=\|\rth(x_{11},x_{12}), \ldots, \rth(x_{n1},x_{n2})\|_{K},
\end{equation}

\noindent
where $\rth$ stands for the counterclockwise rotation of $\R^2$ by the angle $\theta$ with respect to the origin. If a convex body satisfies $(1)$ we will say that \emph{it is invariant with respect to all $\rth$}.

%From the above it is clear that, as in the case of the complex Busemann-Petty problem, here, in order to compare the $(n-1)$-dimensional complex hyperplane sections of two origin symmetric convex bodies we only need to consider those $(2n-2)$-dimensional sections by subspaces that come from complex hyperplanes. Moreover, the convex bodies that appear are only those bodies that when are viewed as compact subsets of $\R^{2n}$ are invariant with respect to any rotation $\rth.$

\medskip

The complex Busemann-Petty problem ([KKZ]) can now be formulated as follows:
Suppose $K$ and $L$ are origin symmetric invariant with respect to all $\rth$
 convex bodies in $\R^{2n}$ such that
 $$\vol_{2n-2}(K\cap H_\xi)\leq \vol_{2n-2}(L\cap H_\xi)$$
 for each $\xi$ from the unit sphere $S^{2n-1}$ of $\R^{2n}.$ Does it follow that
 $$\vol_{2n}(K) \leq \vol_{2n}(L) \  ?$$
 
As it is proved in [KKZ], the answer is affirmative if $n\leq 3$ and negative if $n\geq 4.$ In this article our aim is to extend the result from [KYY] to the complex case.

 \medskip
 
Let $D$ be an origin symmetric convex body in $\C^n.$ For every $\xi \in \C^n, \ |\xi |=1,$ we define the section function 
\begin{equation}\label{eqt:defs}S_{CD}(\xi)=\vol_{2n-2}(D\cap H_{\xi}), \ \forall \xi \in S^{2n-1}.\end{equation}

\noindent
Extending $S_{CD}$ to the whole $\R^{2n}$ as a homogeneous function of degree $-2$ we prove the following:

%\noindent
%Suppose that for two origin symmetric invariant with respect to all $\rth$ convex bodies, $K$ and $L$ in $\R^{2n}$ we have that
%$$\bigl(-\Delta \bigr)^{\alpha/2}S_{CK}(\xi)\leq \bigl(-\Delta \bigr)^{\alpha/2}S_{CL}(\xi), \ \xi \in S^{2n-1},$$
%$\alpha \in \R,$ then the same inequality holds for their $2n$-dimensional volumes.
\begin{main}
Suppose $K$ and $L$ are two origin symmetric invariant with respect to all $\rth$ convex bodies in $\R^{2n}.$ Suppose that $\alpha \in [2n-6, 2n-2), \ n\geq 3.$ If 
\begin{equation}\label{eqt:mSineq}\bigl(-\Delta \bigr)^{\alpha/2}S_{CK}(\xi)\leq \bigl(-\Delta \bigr)^{\alpha/2}S_{CL}(\xi),\end{equation}
for every $\xi \in S^{2n-1}.$ Then

$$\vol_{2n}(K)\leq \vol_{2n}(L).$$
If $\alpha \in (2n-7,2n-6)$ then one can construct two convex bodies $K$ and $L$ that satisfy (\ref{eqt:mSineq}), but $\vol_{2n}(K) > \vol_{2n}(L).$

\end{main}

This means that one needs to differentiate the section functions at least $2n-6$ times and compare the derivatives in order to obtain the same inequality for the volume of the original bodies. Note that if $\alpha=0$ the problem coincides with the original complex Busemann-Petty problem.

\medskip

\section{The Fourier analytic approach}

Throughout this paper we use the Fourier transform of distributions. The Schwartz class of the rapidly decreasing infinitely differentiable functions (test functions) in $\R^n$ is denoted by $\test(\R^n),$ and the space of distributions over $\test(\R^n)$ by $\test^{\prime}(\R^n).$ The Fourier transform $\hat{f}$ of a distribution $f \in \test^{\prime}(\R^n)$ is defined by $\langle\hat{f},\phi\rangle=\langle f,\hat{\phi}\rangle$ for every test function $\phi.$ A distribution is called even homogeneous of degree $p \in \R$ if $\langle f(x),\phi(x/\alpha)\rangle=|\alpha|^{n+p}\langle f,\phi\rangle$ for every test function $\phi$ and every $\alpha \in \R, \ \alpha\neq 0.$ The Fourier transform of an even homogeneous distribution of degree $p$ is an even homogeneous distribution of degree $-n-p.$ A distribution $f$ is called positive definite if, for every test function $\phi, \ \langle f, \phi\ast\overline{\phi(-x)}\rangle \geq 0.$ By Schwartz's generalization of Bochner's theorem, this is equivalent to $\hat{f}$ being a positive distribution in the sense that  $\langle\hat{f},\phi\rangle \geq 0$ for every non-negative test function $\phi.$

We denote by $\Delta$ the Laplace operator on $\R^n,$ and by $|\cdot|_2$ the Euclidean norm in the proper space. Then the fractional powers of the Laplacian are defined by
\begin{equation}\label{eqt:laplacef}
\Bigl(\bigl(-\Delta \bigr)^{\alpha/2}f\Bigr)^{\wedge}=\frac{1}{(2\pi)^n}|x|_2^{\alpha}\hat{f}(x),\end{equation}

\noindent
where the Fourier transform is considered in the sense of distributions. 
%If $\alpha$ is an even integer and $f$ an even distribution, formula (\ref{eqt:laplacef}) gives the negative Laplacian applied $\alpha /2$ times.

 A compact set $K \subset\R^n$ is called a star body, if every straight line that passes through the origin crosses the boundary of the set at exactly two points and the boundary of $K$ is continuous in the sense that the \emph{Minkowski functional} of $K,$ defined by
$$\|x\|_K=\min \{\alpha \geq 0 : x \in \alpha K \}$$
is a continuous function on $\R^n.$
Using polar coordinates it is possible to obtain the following \emph{polar formula of the volume} of the body:

$$\vol_n(K)=\int_{\R^n}\chi (\|x\|_K)dx=\frac{1}{n}\int_{S^{n-1}}\|\theta\|_K^{-n}d\theta.$$

A star body $K$ in $\R^n$ is called $k$-smooth (infinitely sooth) if the restriction of $\|x\|_{K}$ to the sphere $S^{n-1}$ belongs to the class of $C^k(S^{n-1} )\  (C^{\infty}(S^{n-1})).$ It is well-known that one can approximate any convex body in $\R^n$ in the radial metric,
$d(K, L)=\sup \{|\rho_{K}(\xi)-\rho_{L}(\xi)|,\ \xi \in S^{n-1} \},$
by a sequence of infinitely smooth convex bodies. The proof is based on a simple convolution argument (see for example [Sch, Theorem 3.3.1]). It is also easy to see that any convex body in $\R^{2n}$ invariant with respect to all $\rth$ rotations can be approximated in the radial metric by a sequence of infinitely smooth convex bodies invariant with respect to all $\rth.$ This follows from the same convolution argument, because invariance with respect to $\rth$ is preserved under convolutions.
This approximation argument allows us to consider only infinitely smooth origin symmetric convex bodies for the solution to the problem.

If $D$ is an infinitely smooth origin symmetric star body in $\R^n$ and $0<k <n,$ then the Fourier transform of the distribution $\|x\|_D^{-k}$ is a homogeneous function of degree $-n+k$ on $\R^n,$ whose restriction to the sphere is infinitely smooth (see [K5, Lemma 3.16]).

We use a spherical version of Parseval's identity, established in [K3]
(see also [K5, Lemma 3.22]).

 \begin{prop}\label{prop:parseval}
 Let $K$ and $L$ be two infinitely smooth origin symmetric convex bodies in $\R^n$ and $0<p<n.$ Then

 $$\int_{S^{n-1}}\bigl(\|x\|_K^{-p}\bigr)^{\wedge}(\xi)\bigl(\|x\|_L^{-n+p}\bigr)^{\wedge}(\xi)d\xi
 =(2\pi)^n\int_{S^{n-1}}\|x\|_K^{-p}\|x\|_L^{-n+p}dx.$$
\end{prop}

\bigskip

\medskip

Let $H$ be an $(2n-2)-$dimensional subspace of $\R^{2n}$ and $p\leq 2n-2.$  We fix an orthonormal basis $,\{e_1,e_2\},$ in the orthogonal subspace $H^{\perp}.$ For any convex body $D$ in $\R^{2n}$ we define the function $A_{D,H,p}$ as a function on $\R^2$ such that

\begin{equation}\label{def:A}A_{D,H,p}(u)=
\int_{D\cap H_u}|x|_2^{-p} dx, \ u\in \R^2, \end{equation}
where $H_u=\{x\in \R^{2n}:(x,e_1)=u_1,(x,e_2)=u_2\}.$

If the body $D$ is infinitely smooth and $0\leq p <2n-m-2,$ then $A_{D,H,p}$ is $m$-times continuously differentiable near the origin. This can be seen from an argument similar to [K5, Lemma 2.5].

In addition, if we consider the action of the distribution $|u|_2^{-q-2}/ \G(-q/2)$ on $A_{D,H,p}$ we may apply a standard regularization argument (see [GS, p.71-74]) and define the function
\begin{equation}\label{eqt:q-A}
q \longmapsto \left< \frac{|u|_2^{-q-2}}{\G(-\frac{q}{2})},A_{D,H,p}(u)\right>.
\end{equation}

For $q \in \C$ with $ \Re q \leq 2n-p-3,$ the function is an analytic function of $q.$
If $q<0$
\begin{equation}\label{eqt:q<0}\left< \frac{|u|_2^{-q-2}}{\G(-\frac{q}{2})},A_{D,H,p}(u)\right>=
\frac{1}{\G(-q/2)}\int_{\R^2}|u|_2^{-q-2} A_{D,H,p}(u)du.\end{equation}

If $q=2d, \ d \in \N\cup \{0\},$ then
$$\left<\frac{|u|_2^{-q-2}}{\G(-\frac{q}{2})}\Big|_{q=2d}, A_{D,H,p}(u)\right>$$
\begin{equation}\label{eqt:q=2d}=\frac{(-1)^d \pi}{2^{2d}d!}\Delta^d A_{f,D,H}(0),
\end{equation}
where $\Delta=\sum_{i=1}^2 \partial^2/\partial u_i^2$ is the $2$-dimensional Laplace operator (see [GS, p.71-74]).
Note that the function (\ref{eqt:q-A}) is equal, up to a constant, to the fractional power of $\Delta^{q/2}A_{D,H,p}.$ (see [KKZ, p.6-7] or [K4, p.6-7] for complete definition).

\noindent
If the body $D$ is origin symmetric the function $A_{D,H,p}$ is even and for $0< q <2$ we have  (see also [K5, p.39])
$$\left<\frac{|u|_2^{-q-2}}{\G(-\frac{q}{2})}, A_{D,H,p}(u)\right>$$
\begin{equation}\label{eqt:0q2}
=\frac{1}{\G(-\frac{q}{2})}\int_0^{2\pi}
\Bigl(\int_0^{\infty}\frac{A_{D,H,p}(t\theta)-A_{D,H,p}(0)}{t^{1+q}}dt\Bigr)d\theta.
\end{equation}

\medskip

The following proposition is a generalization of [K4], (see also [KKZ, Proposition 4]) with $k=2$. We prove it using a well-known
formula (see for example [GS, p.76]): for any $v\in \R^2$
and $q< -1,$
%$$(v_1^2+v_2^2)^{\frac{-q-2}{2}} $$
\begin{equation} \label{sph}
(v_1^2+v_2^2)^{\frac{-q-2}{2}}= {{\Gamma(-q/2)}\over {2\Gamma((-q-1)/2) \pi^{1/2}}}
\int_{0}^{2\pi} |(v,u)|^{-q-2}\ du.
\end{equation}

\begin{prop}\label{prop:A}
Let $D$ be an infinitely smooth origin symmetric convex body in $\R^{2n}.$ If $-2<q<2n-2, \ 0\leq p \leq 2n-q-3.$ Then for every $(2n-2)$-dimensional subspace $H$ of $\R^{2n}$

$$\left< \frac{|u|_2^{-q-2}}{\G(-\frac{q}{2})}, A_{D,H,p}(u)\right>$$
\begin{equation}\label{eqt:A^{q}}=\frac{2^{-q-2}}{\pi\G\bigl(\frac{q+2}{2}\bigr)(2n-q-p-2)}
\int\limits_{S^{2n-1}\cap H^{\perp}}
\Bigl(\|x\|_D^{-2n+q+p+2}|x|_2^{-p}
\Bigr)^{\wedge}(\theta)d\theta.
\end{equation}
Also, for every $ d \in \N \cup {0}, \ d <n-1$
\begin{equation}\label{eqt:d-A}
\Delta^d A_{D,H,p}(0)=\frac{(-1)^d}{8\pi^2(n-d-1)}\int\limits_{S^{2n-1}\cap H^{\perp}}
\Bigl(\|x\|_D^{-2n+2d+p+2}|x|_2^{-p}
\Bigr)^{\wedge}(\eta)d\eta.
\end{equation}

\end{prop}

\medskip

\textbf{Proof.}
First we assume that $q\in (-2,-1).$ Then 
$$\left< \frac{|u|_2^{-q-2}}{\G(-\frac{q}{2})}, A_{D,H,p}(u)\right>
=\frac{1}{\G(-q/2)}\int_{\R^2}|u|_2^{-q-2} A_{D,H,p}(u)du$$
Using the expression (\ref{def:A}) for the function $A_{D,H,p}$, writing
the integral in polar coordinates and then using (\ref{sph}), we
see that the right-hand side of the latter equation is equal to
$${1\over {\Gamma(\frac{-q}{2})}} \int_{\R^n} \big((x,e_1)^2+(x,e_2)^2)^{\frac{-q-2}{2}}  |x|_2^{-p}\chi(\|x\|_D)\ dx$$
$$={1\over {\Gamma(\frac{-q}{2})(n-q-p-2)}} \int_{S^{n-1}} \left((\theta,e_1)^2+(\theta,e_2)^2 \right)^{\frac{-q-2}{2}} \|\theta\|_D^{-n+q+2}\ d\theta$$
$$={1\over {2\Gamma(\frac{-q-1}{2})\pi^{\frac{1}{2}}(n-q-p-2)}}\ \times$$
$$\int_{S^{n-1}} \|\theta\|_D^{-n+q+p+2} \left(
\int_{0}^{2\pi} \big| ( u_1e_1+u_1e_2, \theta) \big|^{-q-2}\ du \right)\ d\theta$$
$$={1\over {2\Gamma(\frac{-q-1}{2})\pi^{\frac{1}{2}}(n-q-p-2)}}\times$$
\begin{equation} \label{comput}
\int_0^{2\pi} \left(
\int_{S^{n-1}} \|\theta\|_D^{-n+q+p+2} \big| ( u_1e_1+u_2e_2, \theta)
\big|^{-q-2}\ d\theta\right) du.
\end{equation}

Let us show that the function under the integral over $[0,2\pi]$
is the Fourier transform of $\|x\|_D^{-n+q+p+2}|x|_2^{-p}$ at the point
$ u_1e_1+u_2e_2$. For any even test function $\phi\in {\test}(\R^n),$
using the well-known connection between the Fourier and Radon transforms
(see [K5, p.27]) and the
expression for the Fourier transform of the distribution
$|z|_2^{q-1}$ (see [K5, p.38]), we get

$$\langle (\|x\|_D^{-n+q+p+2}|x|_2^{-p})^\wedge, \phi \rangle =
%\langle \|x\|_D^{-n+q+p+2}|x|_2^{-p}, \hat\phi \rangle =
\int_{\R^n} \|x\|_D^{-n+q+p+2}|x|_2^{-p} \hat\phi(x)\ dx $$
$$=\int_{S^{n-1}} \|\theta\|_D^{-n+q+p+2}\left(
\int_0^\infty r^{q+1} \hat\phi(r\theta)\ dr\right) d\theta$$
$$={1\over 2} \int_{S^{n-1}} \|\theta\|_D^{-n+q+p+2}
\Big\langle |r|^{q+1}, \hat\phi(r\theta) \Big\rangle \ d\theta$$
$$={{2^{q+2}\sqrt{\pi}\ \Gamma((q+2)/2)}\over
{2\Gamma((-q-1)/2)}}
\int_{S^{n-1}} \|\theta\|_D^{-n+q+p+2}
\Big\langle |t|^{-q-2},
\int_{(y,\theta)=t} \phi(y)\ dy \Big\rangle \ d\theta$$
$$= {{2^{q+1}\sqrt{\pi} \Gamma((q+2)/2)}\over
{2\Gamma((-q-1)/2)}} \int_{\R^n} \Big(\int_{S^{n-1}}
|(\theta, y)|^{-q-2} \|\theta\|_D^{-n+q+p+2}\ d\theta \Big) \phi(y)\ dy.$$
Since $\phi$ is an arbitrary test function, this
proves that, for every $y\in \R^n\setminus \{0\},$
$$\big(\|x\|_D^{-n+q+p+2}|x|_2^{-p}\big)^\wedge(y)$$
$$ =
{{2^{q+2}\sqrt{\pi} \Gamma((q+2)/2)}\over
{2\Gamma((-q-1)/2)}} \int_{S^{n-1}}
|(\theta, y)|^{-q-2} \|\theta\|_D^{-n+q+p+2}\ d\theta.$$
Together with (\ref{comput}), the latter equality shows that
$$\Big\langle {{|u|_2^{-q-2}}\over {\Gamma(-q/2)}},
A_{D,H,p}(u) \Big\rangle$$
\begin{equation} \label{fracform}
 = {{2^{-q-2}\pi^{-1}}\over {\Gamma((q+2)/2)(n-q-p-2)}}
\int_{S^{n-1}\cap H^\bot} \big(\|x\|_D^{-n+q+p+2}|x|_2^{-p}\big)^\wedge(\theta)
\ d\theta,
\end{equation}
because in our notation $S^{n-1}\cap H^\bot=[0,2\pi].$

We have proved (\ref{fracform}) under the assumption that $q\in (-2,-1).$
However, both sides of (\ref{fracform})
are analytic functions of $q\in \C$ in the domain where
$-2<\Re q<2n-2.$
This implies that the equality (\ref{fracform}) holds for every $q$ from this domain
(see [K5, p.61] for the details of a similar argument).

Putting $q=2m,\ m\in \N\cup \{0\},\ m< n-1$ in (\ref{fracform}) and applying
(\ref{eqt:q=2d}) and the fact that $\Gamma(x+1)=x\Gamma(x)$, we get the second
formula.
\qed

\medskip

\bigskip

Brunn's theorem (see for example [K5, Theorem 2.3]) states that for an origin symmetric convex body and a fixed direction, the central hyperplane section has the maximal volume among all the hyperplane sections perpendicular to the given direction. As a consequence we have the following generalization proved in [KKZ, Lemma 1] for $p=0.$

\begin{prop}\label{prop:brunn}
Suppose $D$ is a $2$-smooth origin symmetric convex body in $\R^{2n},$ then the function $ A_{D,H,p}$ is twice differentiable at the origin and $$\Delta A_{D,H,p}(0) \leq 0.$$
Moreover, for any $q \in (0, 2)$,
$$\left< \frac{|u|_2^{-q-2}}{\G(-\frac{q}{2})}, A_{D,H,p}(u)\right> \geq 0.$$
\end{prop}

\medskip

\textbf{Proof.} Differentiability follows from the same argument as in [K5, Lemma 2.4].

The body $D$ is origin symmetric and convex, so, to prove the first inequality we need to observe that the function $u\longmapsto A_{D,H,p}(u), \ u\in \R^2,$ attains its maximum at the origin:

If $p=0$ then it follows immediately from Brunn's theorem (see [K5, Theorem 2.3] and [KKZ, Lemma 1].)

\noindent
Let $p>0.$ Since $|x|_2^{-p}=p\int_{0}^{\infty}\chi(z|x|_2)z^{q-1}dz,$ we have that for any $u\in \R^2$
\begin{eqnarray*}A_{D,H,p}(u)&=&\int_{D\cap H_u}|x|_2^{-p} dx=p\int_{D\cap H_u}\int_{0}^{\infty}\chi(z|x|_2)z^{q-1}dzdx \\
&=&p\int_{0}^{\infty}z^{q-1}\int_{D\cap H_u}\chi(z|x|_2)dxdz \\
&=&p\int_{0}^{\infty}z^{q-1}\int_{B(1/z)\cap H_u}\chi(\|x\|_D)dx dz,  \end{eqnarray*}

\noindent
where $B(1/z)$ is the unit ball of radius $1/z,.$ Applying Brunn's theorem to the body $B(1/z)\cap D, $ we have that the latter integral is
$$\leq p\int_{0}^{\infty}z^{q-1}\int_{ H}\chi(\|x\|_{B(1/z)\cap D})dx dz=A_{D,H,p}(0).$$

If $q\in (0,2)$ then $ \G(-q/2)<0.$ Hence, for the second inequality we use (\ref{eqt:0q2}) to get that
$$\left<\frac{|u|_2^{-q-2}}{\G(-\frac{q}{2})}, A_{D,H,p}(u)\right>$$
$$
=\frac{1}{\G(-\frac{q}{2})}\int_0^{2\pi}
\Bigl(\int_0^{\infty}\frac{A_{D,H,p}(t\theta)-A_{D,H,p}(0)}{t^{1+q}}dt\Bigr)d\theta \geq 0,$$
since $A_{D,H,p}(u)\leq A_{D,H,p}(0),$ for every $u\in \R^2.$

 \section{Distributions of the form $|x|_2^{-\beta}\|x\|^{-\gamma}$ }\label{distr}

As in the modified real Busemann-Petty problem the solution is closely related to distributions of the form $|x|_2^{-\beta}\|x\|^{-\gamma}.$ 

First, we need a simple observation. The following lemma is crucial for the solution of the problem. 
\medskip

\begin{lemma}\label{modconst}
 For every infinitely smooth origin symmetric invariant with respect to all $\rth$ convex body $D$ in $\R^{2n}$ and every $\xi \in S^{2n-1},$ the Fourier transform of the distribution $|x|_2^{-\beta}\|x\|_{D}^{-\gamma}, \ 0< \beta , \gamma<2n$ is a constant function on $S^{2n-1}\cap H_{\xi}^{\perp}.$
\end{lemma}

\noindent
\textbf{Proof.}
The proof (see [KKZ, Theorem 1], when $\beta=0$) is based on the following observation:

 \noindent
 The body $D$ is invariant with respect to all $\rth.$ So, because of the connection between the Fourier transform and linear transformations, the Fourier transform of $|x|_2^{-\beta}\|x\|_{D}^{-\gamma}$ is also invariant with respect to all $\rth.$ This implies that it is a constant function on $S^{2n-1}\cap H^{\perp}_{\xi}$ because this circle can be represented as the set of all the rotations $\rth,\ \theta \in [0, 2\pi],$ of the vector $\xi \in S^{2n-1}.$

%Let $\xi \in S^{2n-1}.$
As a consequence of the above we have that
\begin{equation}\label{eqt:invariance}
\int_{S^{2n-1}\cap H_{\xi}^{\perp}}\Bigl(|x|_2^{-\beta}\|x\|_{D}^{-\gamma}\Bigr)^{\wedge}(\theta) d\theta=2\pi \Bigl(|x|_2^{-\beta}\|x\|_{D}^{-\alpha}\Bigr)^{\wedge}(\xi).
\end{equation}
\qed

\bigskip

\begin{lemma}\label{lm:posdef}
Let $D$ be an origin symmetric invariant with respect to all $\rth$ convex body in $\R^{2n}, n\geq 3.$ If $q \in (-2,2]$ and $0\leq p< 2n-q-3$ then $|x|_2^{-p}\|x\|_D^{-2n+p+q+2}$ is a positive definite distribution.
\end{lemma}

\medskip

\noindent
\textbf{Proof.}

\medskip

If $p=0$ then by [KKZ, Theorem 3],   $\bigl(\|x\|_D^{-2n+q+2}\bigr)^{\wedge} \geq 0,$ since $2n-q-2 \in [2n-4, 2n).$

\medskip

Let $p>0.$ If $q \in (-2,0)$ then by equation (\ref{eqt:q<0}) and Proposition \ref{prop:A} (formula (\ref{eqt:A^{q}})) we have that
$$\frac{2^{-q-2}}{\pi\G\bigl(\frac{q+2}{2}\bigr)(2n-q-p-2)}
\int\limits_{S^{2n-1}\cap H^{\perp}}
\Bigl(\|x\|_D^{-2n+q+p+2}|x|_2^{-p}
\Bigr)^{\wedge}(\theta)d\theta$$
$$=\frac{1}{\G(-q/2)}\int_{\R^2}|u|_2^{-q-2} A_{D,H,p}(u)du\geq 0.$$
By Lemma \ref{modconst}, the Fourier transform of the distribution 
$|x|_2^{-p}\|x\|_{D}^{-2n+p+q+2}$ is a constant function on $S^{2n-1}\cap H^{\perp}_{\xi}.$ So,

$$\Bigl(|x|_2^{-p}\|x\|_{D}^{-2n+p+q+2}\Bigr)^{\wedge} \geq 0,$$ 
since $\G(\frac{q+2}{2})>0, \ \G(-\frac{q}{2})>0$ and $q<2n-p-2.$

\medskip

Now, if $q=0,$ (\ref{eqt:d-A}) and (\ref{eqt:invariance}) give that
$$A_{D,H,p}(0)=\frac{1}{4\pi (n-1)}\Bigl(|x|_2^{-p}\|x\|_{D}^{-2n+p+q+2}\Bigr)^{\wedge}(\xi) \geq 0.$$

\medskip

For the case where $q \in (0,2)$ we use Proposition \ref{prop:A} and the Remark to get that
$$\left< \frac{|u|_2^{-q-2}}{\G(-\frac{q}{2})}, A_{D,H,p}(u)\right>$$
$$=\frac{2^{-q-1}}{\G\bigl(\frac{q+2}{2}\bigr)(2n-q-p-2)} \Bigl(\|x\|_D^{-2n+q+p+2}|x|_2^{-p}
\Bigr)^{\wedge}(\xi).$$
Then, by the generalization of Brunn's theorem, Proposition \ref{prop:brunn}, the desired follows.

 Lastly, if $q=2,$ (\ref{eqt:d-A}) and (\ref{eqt:invariance}) imply that
 $$\Delta A_{D,H,p}(0)=\frac{-1}{4\pi(n-2)}
\Bigl(\|x\|_D^{-2n+p+4}|x|_2^{-p}
\Bigr)^{\wedge}(\xi).$$
Combining this with Brunn's generalization, since the Laplacian of the function $A_{D,H,p}$ at $0$ is non-positive, we have that
$$\Bigl(|x|_2^{-p}\|x\|_{D}^{-2n+p+4}\Bigr)^{\wedge}(\xi) \geq 0.$$

 \qed

 Before we prove the main result of this article we need the following:

\medskip
\begin{lemma}
 Let $D$ be an infinitely smooth origin symmetric invariant with respect to all $\rth$ convex body in $\R^{2n}$ and $\alpha \in \R.$ Then
 \begin{equation}\label{eqt:laplaces}\bigl(-\Delta \bigr)^{\alpha/2}S_{CD}(\xi)=\frac{1}{4\pi (n-1)}\bigl(|x|_2^{\alpha}\|x\|_{D}^{-2n+2}\bigr)^{\wedge}(\xi)
\end{equation}

\end{lemma}

\noindent
\textbf{Proof.} Let $\xi \in S^{2n-1}.$ As proved in [KKZ, Theorem 1], using the same idea as in Lemma \ref{modconst} (with $r=0$)
\begin{equation}\label{eqt:vol}
\vol_{2n-2}(D\cap H_{\xi})=\frac{1}{4\pi (n-1)}\Bigl(\|x\|_{D}^{-2n+2}\Bigr)^{\wedge}(\xi).
\end{equation}

By the definition of the section function of $D,$ and equation (\ref{eqt:vol}) we obtain the following formula:
\begin{equation}S_{CD}(\xi)=\frac{1}{4\pi (n-1)}\Bigl(\|x\|_{D}^{-2n+2}\Bigr)^{\wedge}(\xi).
\end{equation}

We extend $S_{CD}$ to the whole $\R^{2n}$ as a homogeneous function of degree $-2$ and apply the definition of the fractional powers of the Laplacian. Then, since $\|x\|_{D}^{-2n+2}$ is an even distribution, equation (\ref{eqt:laplaces}) immediately follows.
\qed

\bigskip

                    \section{The solution of the problem.}\label{solMCBP}

\noindent
We consider the affirmative and negative part of the main result separately. The proof follows by the next two theorems.

\begin{thm}\label{thm:pos}(AFFIRMATIVE PART) Let $K$ and $L$ be two infinitely smooth origin symmetric invariant with respect to all $\rth$ convex bodies in $\R^{2n}.$ Suppose that $\alpha \in [2n-6, 2n-2), n\geq 3.$ Then
for every $\xi \in S^{2n-1}$
\begin{equation}\label{eqt:Sineq}\bigl(-\Delta \bigr)^{\alpha/2}S_{CK}(\xi)\leq \bigl(-\Delta \bigr)^{\alpha/2}S_{CL}(\xi)\end{equation}
implies that

$$\vol_{n}(K)\leq \vol_{n}(L).$$
\end{thm}

\noindent
\textbf{Proof.} The bodies $K$ and $L$ are infinitely smooth and invariant with respect to all $\rth$ convex bodies. So by equation (\ref{eqt:laplaces}) the condition in (\ref{eqt:Sineq}) can be written as
\begin{equation}\label{eqt:normineq}\Bigl(|x|_2^{\alpha}\|x\|_{K}^{-2n+2}\Bigr)^{\wedge}\leq \Bigl(|x|_2^{\alpha}\|x\|_{L}^{-2n+2}\Bigr)^{\wedge}.
\end{equation}

 We apply Lemma \ref{lm:posdef} with $p=\alpha$ and $q=2n-\alpha -4$  so that the distribution $|x|_2^{\alpha}\|x\|_{K}^{-2}$ is positive definite. By Bochner's theorem this implies that its Fourier transform is a non-negative function on $\R^{2n}\setminus \{0\}.$ By [K5, Lemma 3.16], it is also continuous, since $K$ is infinitely smooth. Multiply both sides in (\ref{eqt:normineq}) by$\bigl(|x|_2^{-\alpha}\|x\|_{K}^{-2}\bigr)^{\wedge}$ and integrate over the unit sphere $S^{2n-1}.$ Then we can apply Parseval's spherical version, Proposition \ref{prop:parseval}, to get that
 \begin{equation}\label{eqt:volK}
\int_{S^{2n-1}}\|x\|_K^{-2n}dx\leq \int_{S^{2n-1}}\|x\|_K^{-2}\|x\|_L^{-2n+2}.\end{equation}

 \noindent
 Then, by a simple application of H\"older's inequality on formula (\ref{eqt:volK}) and the polar formula of the bodies (see Section 2.) we obtain the affirmative answer to the problem, since
 $$2n\ \vol_{2n}(K)\leq \Bigl(2n\ \vol_{2n}(K)\Bigr)^{1/n}\Bigl(2n\ \vol_{2n}(L)\Bigr)^{(n-1)/n}.$$
\qed

 \medskip

 \bigskip
To prove the negative part we need the following lemma.
 \begin{lemma}\label{lm:nonposdef}
 Let $\alpha \in (2n-7, 2n-6).$ There exists an infinitely smooth origin symmetric convex body $L$ with positive curvature, so that
 $$|x|_2^{-\alpha}\|x\|_{L}^{-2}$$
 is not a positive definite distribution.
 \end{lemma}

We postpone the proof of Lemma \ref{lm:nonposdef} until the end of this section to show that the existence of such a body gives a negative answer to the problem.

\medskip

\begin{thm}\label{thm:neg}(NEGATIVE PART)
Suppose there exists an infinitely smooth, origin symmetric convex body $L$ for which $|x|_2^{-\alpha}\|x\|_{L}^{-2}$ is not a positive definite distribution. Then one can construct an origin symmetric convex body $K$ in $\R^{2n}, n\geq 3,$ so that together with $L$ they satisfy (\ref{eqt:Sineq}), for every $\xi \in S^{2n-1}$ but
$$\vol_{2n}(K)> \vol_{2n}(L).$$
\end{thm}

\noindent
\textbf{Proof.}
The body $L$ is infinitely smooth, so by [K5, Lemma 3.16] the Fourier transform of the distribution $|x|_2^{-\alpha}\|x\|_{L}^{-2}$ is a continuous function on the unit sphere $S^{2n-1}.$ Moreover there exists an open subset $\Omega$ of $S^{2n-1}$ in which $\Bigl(|x|_2^{-\alpha}\|x\|_{L}^{-2}\Bigr)^{\wedge}<0 .$
Since $L$ is invariant with respect to all $\rth$ we may assume that $\Omega$ is also invariant we respect to rotations $\rth$.
 
We use a standard perturbation procedure for convex bodies, see for example [K5, p.96] (similar argument was used in [KKZ, Lemma 5]).
 Consider a non-negative infinitely differentiable even function $g$ supported on $\Omega$ that is also invariant with respect to rotations $\rth.$ We extend it to a homogeneous function of degree $-\alpha-2$ on $\R^{2n}.$ By [K5, Lemma 3.16] its Fourier transform is an even homogeneous function of degree $-2n+\alpha+2$ on $\R^{2n},$ whose restriction to the sphere is infinitely smooth: $\bigl(g(x/|x|_2)|x|_2^{-\alpha-2}\bigr)^{\wedge}(y)=h(y/|y|_2)|y|_2^{-2n+\alpha+2},$ where $h\in C^{\infty}(S^{2n-1}).$

We define a body $K$ so that
$$\|x\|_{K}^{-2n+2}=\|x\|_{L}^{-2n+2}+\varepsilon |x|_2^{-2n+2}h\bigl(\frac{x}{|x|_2}\bigr),$$
for small enough $\varepsilon >0$ so that the body $K$ is strictly convex. Note that $K$ is also invariant with respect to all $\rth.$ We multiply both sides by $\frac{1}{4\pi(n-1)}|x|_2^{\alpha}$ and apply Fourier transform. Then
\begin{eqnarray}\label{eqn:Sc}\bigl(-\Delta \bigr)^{\alpha/2}S_{CK}(\xi)&=&\bigl(-\Delta \bigr)^{\alpha/2}S_{CL}(\xi)+
\frac{\varepsilon(2\pi)^{2n}}{4\pi(n-1)}|x|_2^{-\alpha-2}g\bigl(\frac{x}{|x|_2}\bigr) \\
&\leq & \bigl(-\Delta \bigr)^{\alpha/2}S_{CL}(\xi), \nonumber \end{eqnarray}
since $g$ is non-negative.

On the other hand, we multiply both sides of (\ref{eqn:Sc}) by $\Bigl(|x|_2^{-\alpha}\|x\|_{L}^{-2}\Bigr)^{\wedge}$ and integrate over the sphere,

$$\int_{S^{2n-1}}\Bigl(|x|_2^{-\alpha}\|x\|_{L}^{-2}\Bigr)^{\wedge}(\theta)
\bigl(-\Delta \bigr)^{\alpha/2}S_{CK}(\theta) d\theta$$

\begin{eqnarray*}&=&\int_{S^{2n-1}}\Bigl(|x|_2^{-\alpha}\|x\|_{L}^{-2}\Bigr)^{\wedge}(\theta)
\bigl(-\Delta \bigr)^{\alpha/2}S_{CL}(\theta) d\theta \\
&+& \varepsilon\frac{(2\pi)^{2n}}{4\pi(n-1)}\int_{S^{2n-1}}
\Bigl(|x|_2^{-\alpha}\|x\|_{L}^{-2}\Bigr)^{\wedge}(\theta)g(\theta)d\theta \\
&>& \int_{S^{2n-1}}\Bigl(|x|_2^{-\alpha}\|x\|_{L}^{-2}\Bigr)^{\wedge}(\theta)
\bigl(-\Delta \bigr)^{\alpha/2}S_{CL}(\theta) d\theta,\end{eqnarray*}
since $\Bigl(|x|_2^{-\alpha}\|x\|_{L}^{-2}\Bigr)^{\wedge}<0$ on the support of $g.$
\noindent
Using equation (\ref{eqt:laplaces}) and the spherical version of Parseval's identity, the latter becomes
$$\int_{S^{2n-1}}\|x\|_L^{-2}\|x\|_K^{-2n+2} > \int_{S^{2n-1}}\|x\|_L^{-2n}dx.$$As in Theorem \ref{thm:pos} we apply H\"older's inequality and the polar representation of the volume to obtain the desired inequality for the volumes of the bodies.

\qed

 \noindent
 \textbf{Proof of Lemma \ref{lm:nonposdef}.}
 The construction of the body follows similar steps as in [KYY]. We put $q=2n-\alpha-4,$ so $q \in (2,3).$ From the definition of the fractional derivatives, Proposition \ref{prop:A} and the Remark we see that for a $\xi \in S^{2n-1}$ we need to construct a convex body $D$ so that

$$\int_0^{2\pi}\int_0^{\infty}
t^{-q-1}\Bigl(A_{D,H_{\xi},\alpha}(t\theta)-A_{D,H_{\xi},\alpha}(0)-\Delta A_{D,H_{\xi},\alpha}(0)\frac{t^2}{2}\Bigr)dt d\theta<0$$
since $\G (-\frac{q}{2})>0$ for $q \in (2,3).$

We define the function $$f(|u|)=(1-|u|_2^2-N|u|_2^4)^{\frac{1}{2n-\alpha-2}}, \ u \in \R^2 $$
and consider the body $D$ in $\R^{2n}$ as
$$D=\Big\{(x_{11},x_{12},\ldots ,x_{n1},x_{n2})  \in \R^{2n} : |\bar{x}|_2=|(x_{n1},x_{n2})|_2 \in [-\alpha _{N}, \alpha _{N}],$$
$$\bigl(\sum\limits_{i=1 \atop j=1,2}^{n-1}x^2_{ij}\bigr)^{1/2}\leq f(|\bar{x}|_2)\Big\},
$$
where $a_N$ is the first positive root of the equation $f(t)=0.$
From its definition, the body $D$ is strictly convex with an infinitely smooth boundary. We choose $\xi \in S^{2n-1}$ in the direction of $\bar{x}.$ For $u \in \R^2$ with $|u|_2 \in [0, a_N],$ we write equation (\ref{def:A}) in polar coordinates and get that
\begin{eqnarray*}A_{D,H_{\xi},\alpha}(u)&=&\int_{S^{2n-3}_{u}}
\int_0^{f(|u|_2)}(r^2+|u|_2^2)^{-\frac{\alpha}{2}}r^{2n-3}drd\theta \\
&=&|S^{2n-3}_{u}|\int_0^{f(|u|_2)}(r^2+|u|_2^2)^{-\frac{\alpha}{2}}r^{2n-3}dr,
\end{eqnarray*}
where $|S^{2n-3}_{u}|$ is the volume of the $(2n-3)$-dimensional unit sphere. Note that if $|u|_2 > \alpha_N$ then $A_{D,H_{\xi},\alpha}(u)=0.$ Moreover, if $u=t\theta, \ t\in [0, \infty)$ and $\theta \in S^1,$ the parallel section function $A_{D,H_{\xi},\alpha}(t\theta)$ is independent of $\theta$ since
\begin{equation}\label{eqt:A_t}
A_{D,H_{\xi},\alpha}(t\theta)=|S^{2n-3}_{t}|\int_0^{f(t)}(r^2+t^2)^{-\frac{\alpha}{2}}r^{2n-3}dr.
\end{equation}
Hence, we need to prove that the above construction of the body $D$ gives that
\begin{equation}\label{eqt:countrxample}
\int_0^{\infty}
t^{-q-1}\Bigl(A_{D,H_{\xi},\alpha}(t\theta)-A_{D,H_{\xi},\alpha}(0)-\Delta A_{D,H_{\xi},\alpha}(0)\frac{t^2}{2}\Bigr)dt <0.
\end{equation}
Note that the condition $|u|_2 \in [0, \alpha_N]$ is now equivalent to $t \in [0, \alpha_N].$
In order to prove the above we first compute
$$A_{D,H_{\xi},\alpha}(0)=\frac{|S^{2n-3}|}{2n-\alpha-2}$$
and
$$\Delta A_{D,H_{\xi},\alpha}(0)=-|S^{2n-3}|\Big[\frac{1}{2n-\alpha-2}+\frac{\alpha}{2n-\alpha-4}\Big].$$

Let $\beta_N$ be the positive root of the equation $1-t^2-Nt^4=t^{q+1}.$ We split the integral in
(\ref{eqt:countrxample}) in three parts: $[0, \beta_N], \ [\beta_N,\alpha_N]$ and $[\alpha_N, \infty)$ and work separately. It is not difficult to see that for large $N, \ \alpha_N , \beta_N \simeq N^{-\frac{1}{4}}.$ Also, for every $t \in [0, \alpha_N], \ f(t)>0$ and $f(t)\geq t $ if and only if $t \in [0, \beta_N].$

\medskip

%$A_{D,H_{\xi},\alpha}(t\theta)=|S^{2n-3}_{t}|$

 For the first part, the interval $[0,\beta_N],$ since $f(t)\geq t $, the 2-dimensional parallel section function $A_{D,H_{\xi},\alpha}$ can be easily estimated if we split it into two integrals. For the second we use the inequality $(1+x)^{-\gamma}\leq 1-\gamma x +\frac{\gamma(\gamma+1)}{2}x^2,$ for $\gamma >0$ and $0<x<1.$ Then

 $$\int_0^t(r^2+t^2)^{-\frac{\alpha}{2}}r^{2n-3}dr
 \leq \int_0^tr^{-\alpha+2n-3}dr=\frac{t}{2n-\alpha-2}$$
 and
 $$\int_t^{f(t)}(r^2+t^2)^{-\frac{\alpha}{2}}r^{2n-3}dr\leq
 \int_t^{f(t)} \Bigl[1-\frac{\alpha}{2}\frac{t^2}{r^2}+\frac{\alpha(\alpha+2)}{4}
 \frac{t^4}{r^4}\Bigr]r^{2n-\alpha-3}dr$$

$$=\frac{r^{2n-\alpha-2}}{2n-\alpha-2}-\frac{\alpha}{2}t^2\frac{r^{2n-\alpha-4}}{2n-\alpha-4}
+\frac{\alpha(\alpha+2)}{4}t^4\frac{r^{2n-\alpha-6}}{2n-\alpha-6}\Bigg|_t^{f(t)}$$

$$=\frac{f^{2n-\alpha-2}(t)}{2n-\alpha-2}-\frac{\alpha}{2}t^2\frac{f^{2n-\alpha-4}(t)}{2n-\alpha-4}
+\frac{\alpha(\alpha+2)}{4}t^4\frac{f^{2n-\alpha-6}(t)}{2n-\alpha-6}
-Ct^{2n-\alpha-2},$$

\noindent
where $C=\frac{1}{2n-\alpha-2}-\frac{p}{2(2n-\alpha-4)}+\frac{\alpha(\alpha+2)}{4(2n-\alpha-6)}>0,$ since $n\geq 3$ and $\alpha \in (2n-7,2n-6).$

\noindent
We now use the definition of the function $f$ and the inequality $(1-x)^{\gamma}\geq 1-\gamma x(1-x)^{\gamma-1},$ for $0<\gamma<1$ and $0<x<1.$ We then write

$$=\frac{1-t^2-Nt^4}{2n-\alpha-2}-\frac{\alpha}{2}
\frac{t^2(1-t^2-Nt^4)^{\frac{2n-\alpha-4}{2n-\alpha-2}}}{2n-\alpha-4}$$
$$+\frac{\alpha(\alpha+2)}{4}\frac{t^4(1-t^2-Nt^4)^{\frac{2n-\alpha-6}{2n-\alpha-2}}}
{2n-\alpha-6}
-Ct^{2n-\alpha-2}$$
$$\leq \frac{1-t^2-Nt^4}{2n-\alpha-2}-\frac{\alpha t^2}{2(2n-\alpha-4)}+
\frac{\alpha t^2(t^2+Nt^4)}{2(2n-\alpha-2)(1-t^2-Nt^4)^{\frac{2}{2n-\alpha-2}}}$$
$$+
\frac{\alpha(\alpha+2)t^4}{4(2n-\alpha-6)}-
\frac{\alpha(\alpha+2)}{4(2n-\alpha-2)}\frac{t^4(t^2+Nt^4)}{(1-t^2-Nt^4)^{\frac{4}{2n-\alpha-2}}}-Ct^{2n-\alpha-2}.$$

Hence, we have that
$$\int_0^{\beta_N}
t^{-q-1}\Bigl(A_{D,H_{\xi},\alpha}(t\theta)-A_{D,H_{\xi},\alpha}(0)-\Delta A_{D,H_{\xi},\alpha}(0)\frac{t^2}{2}\Bigr)dt$$

$$=\int_0^{\beta_N}
t^{-q-1}\Bigl(Ct^{2n-\alpha-2}-Dt^4+E\frac{t^2(t^2+Nt^4)}{(1-t^2-Nt^4)^{\frac{2}{2n-\alpha-2}}}$$
\begin{equation}\label{betaN}
-F\frac{t^4(t^2+Nt^4)}{(1-t^2-Nt^4)^{\frac{4}{2n-\alpha-2}}}\Bigr)dt, \end{equation}

where  $E=\frac{\alpha}{2(2n-\alpha-2)}>0,$
$F=\frac{\alpha(\alpha+2)}{2n-\alpha-2}>0$ and $D=\frac{N}{2n-\alpha-2}-\frac{\alpha(\alpha+2)}{4(2n-\alpha-6)}>0, $ for $N$ large enough.

Now, in order to obtain an upper bound for (\ref{betaN}) we need to estimate four different integrals. The first one simply gives $\frac{C}{2}\beta_N^2\simeq C_1 N^{-\frac{1}{2}},$ and the second  $\frac{D}{4-q}\beta_N^{-q+4}\simeq D_1 N^{\frac{q-4}{4}},$ for large $N.$ For the third one, we make a change of variables, $u=N^{\frac{1}{4}}t$ and get
$$E\int_0^{\beta_N}\frac{t^{-q+1}(t^2+Nt^4)}{(1-t^2-Nt^4)^{\frac{2}{2n-\alpha-2}}}dt=
EN^{\frac{q-2}{4}}\int\limits_0^{\beta_N N^{\frac{1}{4}}}
\frac{u^{-q+3}(N^{-\frac{1}{2}}+u^2)}{(1-u^2N^{-\frac{1}{2}}-u^4)^{\frac{2}{2n-\alpha-2}}}du$$
$$\leq E_1N^{\frac{q-2}{4}},$$
since $\beta_N N^{\frac{1}{4}}\longrightarrow 1$ as $N \rightarrow \infty$ and the integral $\int_0^1
\frac{u^{-q+5}}{(1-u^4)^{\frac{2}{2n-\alpha-2}}}du$ converges.

We apply the same change of variables, $u=N^{\frac{1}{4}}t,$ for the last integral and find that it is comparable to $ N^{\frac{q}{4}-1}.$  
\begin{equation}\label{F}F\int_0^{\beta_N}\frac{t^{-q+3}(t^2+Nt^4)}{(1-t^2-Nt^4)^{\frac{4}{2n-\alpha-2}}}dt=
FN^{\frac{q}{4}-1}\int\limits_0^{\beta_N N^{\frac{1}{4}}}
\frac{u^{-q+3}(u^2N^{-\frac{1}{2}}+u^4)}{(1-u^2N^{-\frac{1}{2}}-u^4)^{\frac{4}{2n-\alpha-2}}}du.
\end{equation}
The integrand in the latter is a positive increasing function of $u$ and $\beta_N N^{\frac{1}{4}}\longrightarrow 1$ as $N \rightarrow \infty.$ So, we can roughly bound the integral from below by a positive constant and have that equation (\ref{F}) is greater than $F_1N^{\frac{q}{4}-1},$ where $F_1>0.$
%$$FN^{\frac{q}{4}-1}\int\limits_0^{\beta_N N^{\frac{1}{4}}}
%\frac{u^{-q+3}(u^2N^{-\frac{1}{2}}+u^4)}{(1-u^2N^{-\frac{1}{2}}-u^4)^{\frac{4}{2n-\alpha-2}}}du>FF_1$$
%Then, as before, $\beta_N N^{\frac{1}{4}}\longrightarrow 1$ as $N \rightarrow \infty$ and the integral $\int_0^1
%\frac{u^{-q+7}}{(1-u^4)^{\frac{4}{2n-\alpha-2}}}du$ converges.

\medskip

 In the second interval, we use the fact that $A_{D,H,\alpha}(t\theta) \leq A_{D,H,\alpha}(0)$ since central sections have maximum volume. Then, since $t<<1,$ we have that
$$\int_{\beta_N}^{\alpha_N}
t^{-q-1}\Bigl(A_{D,H_{\xi},\alpha}(t\theta)-A_{D,H_{\xi},\alpha}(0)-\Delta A_{D,H_{\xi},\alpha}(0)\frac{t^2}{2}\Bigr)dt$$
$$\leq \int_{\beta_N}^{\alpha_N}
t^{-q-1}\Bigl(\frac{1}{2n-\alpha-2}+\frac{\alpha}{2(2n-\alpha-4)}\Bigr)t^2dt
<A\int_{\beta_N}^{\alpha_N}t^{-q-1}dt.$$

\noindent
Recall that $\alpha_N$ and $\beta_N$ are the positive solutions of the equations $f(t)=0$ and
$1-t^2-Nt^4=t^{q+1}$ respectively, and that for large $N, \ \alpha_N\simeq N^{-\frac{1}{4}}.$ Then, it is
not difficult to see that
$$A\int_{\beta_N}^{\alpha_N}t^{-q-1}dt \leq \frac{A}{(\alpha_N+\beta_N)(1+N(\alpha_N^2+\beta_N^2))}
\simeq AN^{-\frac{1}{4}},$$

\noindent
see [KYY, p.204] for details.

\medskip

 Lastly, for the interval $[\alpha_N, \infty),$ we use the fact that $A_{D,H_{\xi},\alpha}(t\theta)=0.$ Then, we have that
$$\int_{\alpha_N}^{\infty}t^{-q-1}\Bigl(-A_{D,H_{\xi},\alpha}(0)-\Delta A_{D,H_{\xi},\alpha}(0)\frac{t^2}{2}\Bigr)dt$$
$$=\int_{\alpha_N}^{\infty}\Bigl[-\frac{t^{-q-1}}{2n-\alpha-2}+\Bigl(\frac{2}{2n-\alpha-2}+\frac{\alpha}{2n-\alpha-4}\Bigr)\frac{t^{-q+1}}{2}\Bigr]dt$$
$$=-A_1\alpha_N^{-q}+A_2\alpha^{-q+2}\simeq -A_1N^{\frac{q}{4}}+A_2N^{\frac{q-2}{4}},$$

\noindent
where $A_1,A_2 >0.$

Combining all the above estimations, for $N$ large enough, we obtain the following upper bound for the integral in (\ref{eqt:countrxample}), 

$$\int_0^{\infty}
t^{-q-1}\Bigl(A_{D,H_{\xi},\alpha}(t\theta)-A_{D,H_{\xi},\alpha}(0)-\Delta A_{D,H_{\xi},\alpha}(0)\frac{t^2}{2}\Bigr)dt d\theta$$
$$< C_1 N^{-\frac{1}{2}}+D_1 N^{\frac{q-4}{4}}+E_1N^{\frac{q-2}{4}}-F_1 N^{\frac{q}{4}-1}+AN^{-\frac{1}{4}}
-A_1N^{\frac{q}{4}}+A_2N^{\frac{q-2}{4}},$$
which clearly shows that it is negative since all the constants are positive and $q\in (2,3).$

\qed

\bigbreak

{\bf Acknowledgments:}
The author was partially supported by the
NSF grant DMS-0652571. Also, part of this work was carried out when the author was visiting the Department of Mathematics of University of Crete, which the author thanks
for its hospitality.

 \end{document}